\RequirePackage{fix-cm}
\documentclass[smallextended]{svjour3}       
\smartqed  
\usepackage{graphicx}
\usepackage{comment}
\usepackage{amsfonts,amsmath,amssymb,graphicx,tabularx,booktabs}
\usepackage{mathtools}

\newcommand{\bb}{\begin{equation}}
\newcommand{\ee}{\end{equation}}
\newcommand{\Z}{{\mathbb Z}}

\newcommand{\ICG}{\mathrm{ICG}}

\newcommand{\maxspread}{\mathrm{maxspread}}

 \newtheorem{thm}{Theorem}
 
 \newtheorem{lem}[thm]{Lemma}

\newcommand{\QED} {\hfill$\square$}

\DeclareMathOperator{\lcm}{lcm}

\newenvironment{dok} {\par \noindent \textbf{Proof. }}{\QED \par \bigskip \par}

\begin{document}

\title{The least eigenvalues of integral circulant graphs \thanks{The authors gratefully acknowledge
support from the research project of the Ministry of Education, Science and Technological Development of the Republic of Serbia (number 451-03-47/2023-01/ 200124).
}
}


\author{Milan Ba\v si\'c        
}


\institute{M. Ba\v si\'c \at Faculty of Sciences and Mathematics, University of Ni\v{s}
              Vi\v segradska 33, 18000 Ni\v s, Serbia \\
              \email{basic\_milan@yahoo.com}           
}

\date{Received: date / Accepted: date}

\maketitle

\begin{abstract}
The integral circulant graph $\ICG_n (D)$ has the vertex set \\$Z_n = \{0, 1, 2, \ldots, n - 1\}$, where vertices $a$ and $b$ are adjacent if $\gcd(a-b,n)\in D$, with $D \subseteq \{d : d \mid n,\ 1\leq d<n\}$.
In this paper, we establish that the minimal value of the least eigenvalues (minimal least eigenvalue) of integral circulant graphs $\ICG_n(D)$, given an order $n$ with its prime factorization $p_1^{\alpha_1}\cdots p_k^{\alpha_k}$, is equal to $-\frac{n}{p_1}$.
Moreover, we show that the minimal least eigenvalue of connected integral circulant graphs $\ICG_n(D)$ of order $n$ whose complements are also connected is equal to $-\frac{n}{p_1}+p_1^{\alpha_1-1}$.
Finally, we determine the second minimal eigenvalue among all least eigenvalues within the class of connected integral circulant graphs of a prescribed order $n$ and show it to be equal to $-\frac{n}{p_1}+p_1-1$ or $-\frac{n}{p_1}+1$, depending on whether $\alpha_1>1$ or not, respectively. In all the aforementioned tasks, we provide a complete characterization of graphs whose spectra contain these determined minimal least eigenvalues.


\end{abstract}

\section{Introduction}
\label{intro}
Circulant graphs are Cayley graphs over a cyclic group. A graph is
called integral when all the eigenvalues found within its adjacency matrix are integers.
Integral graphs have been the subject of extensive research within the literature, 
with considerable effort directed towards identifying specific classes of graphs possessing integral spectra. 
Over the past two decades, circulant graphs have garnered increasing interest within the field of graph theory and its practical applications. 
Their presence extends into diverse areas such as coding theory, VLSI design, Ramsey theory, and numerous other domains.
Due to their possession of numerous noteworthy characteristics, including vertex transitivity referred to as mirror symmetry, 
circulant graphs find practical application in the field of quantum information transmission.
They are also proposed as viable models for quantum spin networks, which facilitate the phenomenon of perfect state transfer, as documented in \cite{Ba14}.
In the context of quantum communication, a key attribute of these graphs, particularly those exhibiting an integral spectrum, resides in their capability to faithfully transmit quantum states while preserving the network's topological integrity. Notably, integral circulant graphs and unitary Cayley graphs, which constitute a subset of integral circulants, have discovered meaningful applications in the field of molecular chemistry, particularly in the modeling of energy-like quantities  such as the heat of formation of hydrocarbons, as reported in references \cite{BlSh08,RaVe09}.

In this paper we study connected graphs whose least
eigenvalues are minimal and second minimal among all integral circulant graphs of a prescribed order. 
These particular eigenvalues are referred to as the {\it minimal least eigenvalues} and {\it second minimal least eigenvalues}.
In the existing body of literature, there are findings which are related to the minimal least eigenvalue problem, where the objective is to identify connected graphs within a specified class of graphs of a prescribed order that possess the lowest possible least eigenvalue.
For instance, in \cite{FanWangGao08}, the authors find the unique graph possessing the minimal least eigenvalue among the class of connected unicyclic graphs having an order of $n$ and girth $k$.
Additionally, they examine the least eigenvalues of graphs within the set of connected unicyclic graphs of a specified order $n$.
Furthermore, the unique graphs exhibiting the minimal least eigenvalue among all connected bicyclic graphs, cacti, 
and cacti with pendent vertices of order $n$ are ascertained in separate studies \cite{PeBoAl09,PeAlSi11,XingZhou13}.
Nevertheless, the quest for determining the minimal least eigenvalue remains compelling when the scope is restricted to specific classes of connected graphs with a given order, such that their complement graphs also maintain connected.
The authors in \cite{Wang2019} systematically ascertain and categorize all graphs belonging to the class of connected graphs of a given order, the complements of which exhibit cacti-like properties. 
Furthermore, the characterization of all connected graphs of a predetermined order, such that their complements are connected unicyclic graphs, is provided in the publication referenced as \cite{WangFanLiZhang15}.


The graph's spread is defined as the maximum modulo difference observed among all pairs of eigenvalues originating from the graph's adjacency matrix.
The task of determining graphs with the maximum spread within a specific class of connected graphs of a prescribed order is referred to as the {\it maximum spread problem}.
This maximum spread problem has a close relationship with both the maximum index problem and the minimal least eigenvalue problem within a given set of graphs.
In the papers referenced as \cite{FanWangGao08,PeBoAl09,AlPe2015}, the authors determine the unique graphs that possess the maximum spread among all connected unicyclic graphs, bicyclic graphs, and cacti of order $n$.


In this paper, we undertake an examination of the spectral properties of (integral) circulant graphs, a subject that has been the focus of numerous recent studies. 
For instance, characterizing graphs with maximum spread, determining graphs with a limited number of eigenvalues, investigating the multiplicativity of specific eigenvalues, 
analyzing the number of spanning trees, and studying cospectrality, among others, have been explored in \cite{BaIl15,Ba22,Ba23,Ba24,Sander2018,SanderSander2015,SanderSander2023}.
However, recently, there has been notable attention within the research community towards the structural properties of integral circulant graphs as well.
For example, maximal diameter, geometric kernel, clique polynomial, chromatic polynomial, distance polynomial, etc., have been the subjects of inquiry in \cite{BaIl15,BaIlSt23,Sander2021}.

Following an introductory section where we establish the necessary notation and concepts pertaining to integral circulant graphs, the third section presents Lemma 1, which provides a lower bound for the sum of Ramanujan functions, which are dependent on a specific divisor set.
This crucial fact provides the foundation for deriving a tight lower bound applicable to any eigenvalue associated with an integral circulant graph, as expounded in Theorem 2. 
Notably, this lower bound represents the minimal least eigenvalue within the class of connected integral circulant graphs of a prescribed order. 
The associated theorem also give a characterization of the unique graph whose adjacency matrix spectrum encompasses this minimal least eigenvalue.
Additionally, subsequent to establishing this assertion, we undertake an investigation into the maximum spread of integral circulant graphs, using the minimal least eigenvalue. 
In Theorem 3, we provide a comprehensive characterization of integral circulant graphs that exhibit the maximum spread.
The characterization of graphs featuring the minimal least eigenvalue among all integral circulant graphs $\ICG_n(D)$ of a given order $n$, with the additional condition that both $\ICG_n(D)$ and its complement $\overline{\ICG_n(D)}$ are connected, is presented in Theorem 4. Lastly, we determine the second minimal among all least eigenvalues within the class of connected integral circulant graphs of a prescribed order $n$ and enumerate all such graphs whose spectra include this specific eigenvalue (Theorem 5).
The proofs necessitate a thorough and exhaustive discussion, encompassing a multitude of distinct cases.


\section{Preliminaries}
\label{sec:1}

 A {\it circulant graph} $G(n;S)$ is a
graph on vertices $\Z_n=\{0,1,\ldots,n-1\}$ such that vertices $i$
and $j$ are adjacent if and only if $i-j \equiv s \pmod n$ for some
$s \in S$. A set $S$ is called the {\it symbol} of graph $G(n;S)$.
In our analysis, we specifically consider undirected graphs that lack self-loops. Thus, we make the assumption that $S$ is equal to its complement with respect to $n$, denoted as $n-S$, such that $S=n-S=\{n-s\ |\ s\in S\}$, and additionally, we assume that $0$ does not belong to set $S$. It is worth emphasizing that the degree of the graph $G(n;S)$ is equal to the cardinality of set $S$, denoted as $|S|$.

The eigenvalues and corresponding eigenvectors of the graph $G(n;S)$ can be expressed as follows

\begin{equation*} \label{eq:eigenvalues unwigted} \lambda_j=\sum_{s \in S}
\omega^{js}_n, \quad v_j=[1 \ \omega_n^j \ \omega_n^{2j} \cdots
\omega_n^{(n-1)j}]^T, 
\end{equation*}
 where
$\omega_n=e^{i\frac{2\pi}n}$ is the $n$-th root of unity
\cite{Davis70}.

\smallskip

Circulant graphs represent a specific subset within the broader class of Cayley graphs.
Consider a multiplicative group denoted as $\Gamma$ with an identity element represented as $e$.
For a subset, denoted as $S$, of $\Gamma$ where $e\not\in S$, and where the set $S$ is its own inverse, expressed as $S^{-1} = \{s^{-1}\ |\ s\in S\}=S$, the Cayley graph $X = Cay(\Gamma,S)$ is defined as an undirected graph. This graph has a vertex set denoted as $V(X)=\Gamma$ and an edge set denoted as $E(X) = \{{a,b}\ |\ ab^{-1}\in S\}$. 
It can be observed that a graph qualifies as circulant when it is a Cayley graph associated with a cyclic group, implying that its adjacency matrix is cyclic.

A graph is considered {\it integral} when all of its eigenvalues are integers. 
A circulant graph $G(n;S)$ is integral if and only if
$$
S=\bigcup_{d \in D} G_n(d),
$$
for some set of divisors $D \subseteq D_n$ \cite {wasin}. Here
$G_n(d)=\{ k \ : \ \gcd(k,n)=d, \ 1\leq k \leq n-1 \}$, and $D_n$ is
the set of all divisors of $n$, different from $n$.

Therefore an {\it integral circulant graph}
$G(n;S)$ is
defined by its order $n$ and the set of divisors $D$. 
An integral circulant graph with $n$ vertices(in further text $\ICG$), defined by the set of
divisors $D \subseteq D_n$ will be denoted by $\ICG_n(D)$.
Unitary Cayley graph is defined as $\ICG_n(\{1\})$.
From the
above characterization of integral circulant graphs we have that the
degree of an integral circulant graph is $\deg \ICG_n(D)=\sum_{d \in
D}\varphi(\frac{n}{d}). $ Here $\varphi(n)$ denotes the Euler-phi function
\cite{HardyWright}. If $D=\{d_1,\ldots,d_k\}$, it is easy to see
that $\ICG_n(D)$ is connected if and only if
$\gcd(d_1,\ldots,d_k)=1$, given that $G(n;s)$ is connected if and
only if $\gcd(n, S)=1$\ \cite{Hwang03}.

Let us denote the expression as $c(j,n)$, defined by the following equation
\begin{equation}
c(j,n)=\mu(t_{n,j})\frac{\varphi(n)}{\varphi(t_{n,j})}, \quad
t_{n,j}=\frac n{\gcd(n,j)}.
 \label{ramanujan}
\end{equation} 

Here, the function $\mu$ referred to is the Möbius function, which is defined as follows

\begin{eqnarray*}
\mu(n)&=&\left\{
\begin{array}{rl}
1, &  \mbox{if}\ n=1  \\
0, & \mbox{if $n$ is not square--free} \\
(-1)^k, & \mbox {if $n$ is product of $k$ distinct prime numbers}.
\end{array} \right.
\end{eqnarray*}
The expression $c(j,n)$ is known as the {\it Ramanujan function}
(\cite[p.~55]{HardyWright}). From (\ref{ramanujan}) it can be easily seen that $c(j,n)=c(n-j,n)$, for $0\leq j\leq n$.   The eigenvalues of $\ICG_n(D)$  can be
expressed in terms of the Ramanujan function as follows
\begin{equation}
\lambda_j(n,D)=\sum_{d\in D} c(j,\frac{n}{d}).
\label{ldef} 
\end{equation}

When it becomes evident within the context that the eigenvalues $\lambda_j(n, D)$ pertain to the graph $\ICG_n(D)$, we will exclude the parameters $n$ and $D$ from the notation, and simply denote it as $\lambda_j$.
%


%

In this paper, we establish that the order of $\ICG_n(D)$ can be expressed as the prime factorization $n=p_1^{\alpha_1}\cdots
p_k^{\alpha_k}$. Furthermore, when considering a prime number $p$ and an integer $n$, we define $S_p(n)$ as the maximum exponent $\alpha$ for which $p^{\alpha} \mid n$.

\section{Minimal least eigenvalues of integral circulant graphs}

For an arbitrary divisor  $d_0$ of $n$ such that
$p_1\nmid d_0$ we denote $D_n^{d_0}=\{p_1^{\beta_1}d_0\ |\ 0\leq
\beta_1\leq \alpha_1\}$.

\begin{lem}
\label{lem:mineig}
For an arbitrary subset $D\subseteq D_n^{d_0}$ and $0\leq j\leq n-1$
holds that
$$
\sum_{d\in D} c(j,\frac{n}{d})\geq -\frac{\varphi(\frac{n}{d_0})}{p_1-1}.
$$
\end{lem}
\begin{dok}
Let $1\leq j\leq \lfloor n/2 \rfloor$  be an arbitrary index such
that $j=p_1^{\gamma_1}j_1$, for $p_1\nmid j_1$.

Without loss of generality, we may assume that $0\leq \gamma_1\leq
\alpha_1$. Indeed, if $\gamma_1>\alpha_1$ then
$n-j=p_1^{\alpha_1}(n_1-p_1^{\gamma_1-\alpha_1}j_1)$, where we define $n=p_1^{\alpha_1}n_1$. Therefore, in that case we
can carry out the proof using the index $n-j$ instead of $j$. This is due to the properties that $c(j,\frac{n}{d})=c(n-j,\frac{n}{d})$ and $S_{p_1}(n-j)=\alpha_1$.

Now, for an arbitrary $d\in D$ such that $d=p_1^{\beta_1}d_0$ and
$0\leq \beta_1\leq \alpha_1$ we have that
\begin{eqnarray}
\label{eq: t} t_{\frac{n}{d},j}&=& \frac{p_1^{\alpha_1 }n_1}{p_1^{\beta_1}d_0
\gcd(p_1^{\alpha_1-\beta_1}n_1/d_0,p_1^{\gamma_1}j_1)}\nonumber\\
&=& p_1^{\alpha_1-\beta_1-\min(\alpha_1-\beta_1,\gamma_1)}\frac{n_1}{d_0\gcd(n_1/d_0,j_1)}.
\end{eqnarray}

In the rest of the proof, t simplify, let us use the notation $N = N(n_1, d_0, j_1) = t_{\frac{n_1}{d_0},j_1} = \frac{n_1}{d_0\gcd(\frac{n_1}{d_0},j_1)}$. Additionally, it is important to observe that $p_1$ does not divide $N$.

Now, we distinguish three cases depending on the values of $\beta_1$.

\medskip

{\bf Case 1.} $\beta_1 \leq \alpha_1-\gamma_1-2$. It directly implies
that $\gamma_1 \leq \alpha_1-\beta_1-2$ and $\alpha_1-\beta_1-\gamma_1\geq 2$.
Consequently, based on equation (\ref{eq: t}), we can deduce that $p_1^2$ divides $t_{\frac{n}{d},j}$ since $p_1^{\alpha_1 - \beta_1 - \gamma_1}$ divides $t_{\frac{n}{d},j}$.
As a result, according to (\ref{ramanujan}), $c(j, \frac{n}{d}) = 0$.

\smallskip

{\bf Case 2.} $\beta_1 = \alpha_1-\gamma_1-1$. Since $\gamma_1 =
\alpha_1-\beta_1-1$, from (\ref{eq: t}) it holds that $t_{\frac{n}{d},j}= p_1N$  and
hence, following (\ref{ramanujan}), 
$$
c(j,\frac{n}{d}) =
\mu(p_1N)\frac{\varphi(p_1^{\alpha_1-\beta_1}n_1/d_0)}{(p_1-1)\varphi(N)}=-\mu(N)\frac{p_1^{\gamma_1}\varphi(\frac{n_1}{d_0})}{\varphi(N)}.
$$

\smallskip

{\bf Case 3.} $\alpha_1-\gamma_1\leq \beta_1\leq\alpha_1$. The
previous assumption is equivalent to $\alpha_1-\beta_1\leq
\gamma_1\leq\alpha_1$, leading us to deduce that $t_{\frac{n}{d},j}=N$ and

\begin{eqnarray}
\label{eq:Case 3}
c(j,\frac{n}{d})&=&
\mu(N)\frac{\varphi(p_1^{\alpha_1-\beta_1}n_1/d_0)}{\varphi(N)}\nonumber\\
&=&\left\{
\begin{array}{rl}
\mu(N)\frac{\varphi(\frac{n_1}{d_0})}{\varphi(N)}, &  \alpha_1 =\beta_1 \\
\mu(N)\frac{p_1^{\alpha_1-\beta_1-1}(p_1-1)\varphi(\frac{n_1}{d_0})}{\varphi(N)}, & \alpha_1-\gamma_1\leq\beta_1<\alpha_1 \\
\end{array}\right..
\end{eqnarray}

For a given $j$, we see that for all divisors $d\in D$ that satisfy
the condition from the Case 3, the values of the function $c(j,\frac{n}{d})$ have the same signs as $\mu(N)$. On the other hand, if $N$ is a square
free number, $c(j,\frac{n}{d})$ and $\mu(N)$ have the opposite signs for the
divisors $d\in D$ that satisfy the condition from Case 2.

\smallskip

If $\mu(N)=-1$, from the discussion above, we conclude that the minimal value of $\sum_{d\in D} c(j,\frac{n}{d})$ can be obtained for all divisors $d$ satisfying the condition in Case 3:
\begin{eqnarray}
 \sum_{d\in D} c(j,\frac{n}{d})&\geq&
\sum_{\beta_1=\alpha_1-\gamma_1}^{\alpha_1}
c(j,\frac{n}{p_1^{\beta_1}d_0})\nonumber\\
&=&\mu(N)\frac{\varphi(\frac{n_1}{d_0})}{\varphi(N)}(1+(p_1-1)(1+p_1+\cdots+p_1^{\gamma_1-1})) \nonumber \\
&=&-\frac{\varphi(\frac{n_1}{d_0})}{\varphi(N)}p_1^{\gamma_1}\geq-\frac{\varphi(\frac{n_1}{d_0})}{\varphi(N)}p_1^{\alpha_1}=-
\frac{\varphi(\frac{n}{d_0})p_1}{\varphi(N)(p_1-1)}\geq -
\frac{\varphi(\frac{n}{d_0})}{p_1-1}.
\label{eq:sum c}
\end{eqnarray}

\smallskip

The final inequality is valid since $p_1\leq
\varphi(N)$. This is evident because of $\mu(N)=-1$ and $p_1\nmid N$, which lead us to the conclusion that
$N\geq p_2$ and, consequently, $p_1\leq \varphi(p_2)\leq
\varphi(N)$.
This implies that the equality in (\ref{eq:sum c}) is satisfied if and only if $\gamma_1=\alpha_1$, $p_1=2$, and $N=p_2=3$. 
Moreover, from $\gamma_1=\alpha_1$, we have that $0\leq \beta_1\leq \alpha_1$, which is equivalent to $D_n^{d_0}=D$.

%

\medskip

If $\mu(N)= 1$, we deduce from Case 2 
\begin{equation}
\label{eq:sumc1}
\sum_{d\in D} c(j,\frac{n}{d})\geq
-\frac{p_1^{\alpha_1-\beta_1-1}\varphi(\frac{n_1}{d_0})}{\varphi(N)}\geq
-p_1^{\alpha_1-1}\varphi(\frac{n_1}{d_0})= -\frac{\varphi(\frac{n}{d_0})}{p_1-1}.
\end{equation}
Observe that the equality is satisfied exclusively when $\beta_1=0$ and $N=1$, which is the case when $\gamma_1=\alpha_1-\beta_1-1=\alpha_1-1$.

\medskip

Notice that the case $\mu(N) = 0$ has no impact on the inequality intended to be shown.

\end{dok}

By $\overline{D}_{p_1}$ we denote all divisors $d\in D_n$ such that
$p_1\nmid d$, i.e. $\overline{D}_{p_1}=\{d\in D_n\ |\ p_1\nmid d\}$.

\begin{thm}
\label{thm:mleigenvalues}
Let $\ICG_n(D)$ be an arbitrary integral circulant graph of order $n$. For every $1\leq j\leq n-1$, the following inequality is valid
$$
\lambda_j(n,D)\geq -\frac{n} {p_1},
$$ 
with equality occurring if and only if $D=\overline{D}_{p_1}$ and $\frac{n}{p_1}\mid j$.
\end{thm}
\begin{dok}
Notice that we can represent the divisor set $D$ in a more convenient form as follows
$$
D=\bigcup_{ d_0\in \overline{D}_{p_1}\cap D} (D_n^{d_0}\cap D).
$$
Observe that the sets $D_n^{d_0}\cap D$, for every $d_0\in \overline{D}_{p_1}\cap D$, are pairwise disjoint. 
If we denote $n_1=\frac{n}{p_1^{\alpha_1}}$,
according to Lemma \ref{lem:mineig}, it can be concluded that
\begin{eqnarray}
\label{ineq:lambdaj} \lambda_j(n,D)&=&
\sum_{d_0\in\overline{D}_{p_1}\cap D} \sum_{d\in D_n^{d_0}\cap D}
c(j,\frac{n}{d})\geq -\sum_{d_0\in\overline{D}_{p_1}}
\frac{\varphi(\frac{n}{d_0})}{p_1-1}\nonumber \\
&=&-\sum_{d_0\in\overline{D}_{p_1}}
\frac{p_1^{\alpha_1-1}(p_1-1)\varphi(\frac{n_1}{d_0})}{p_1-1}=-p_1^{\alpha_1-1}n_1=-\frac
{n} {p_1}.
\end{eqnarray}

The equality in 
(\ref{ineq:lambdaj}) can be  attained if and only if
$\overline{D}_{p_1}\subseteq D$ and \\$\sum_{d\in D_n^{d_0}\cap D} c(j,\frac{n}{d}) =
-\frac{\varphi(\frac{n}{d_0})}{p_1-1}$, for every $d_0\in \overline{D}_{p_1}$. Furthermore, it is possible that the equality $\sum_{d\in D_n^{d_0}\cap D} c(j,\frac{n}{d}) =
-\frac{\varphi(\frac{n_1}{d_0})}{p_1-1}$ is achieved for the divisors $d\in D_n^{d_0}\cap D$ that satisfy either Cases 2 or Case 3 of Lemma \ref{lem:mineig}.

From Case 3 of the proof of Lemma \ref{lem:mineig} we see that the
equality \\$\sum_{d\in D_n^{d_0}\cap D} c(j,\frac{n}{d}) =
-\frac{\varphi(\frac{n}{d_0})}{p_1-1}$ holds, for every $d_0\in \overline{D}_{p_1}$ if and only if $p_1=2$,
$\frac{n_1}{d_0\gcd(n_1/d_0,j_1)}=3$ and $D_{n}^{d_0}\subseteq D$, where $j_1=\frac{j}{p_1^{\gamma_1}}$.
This means that $\overline{D}_{p_1}\cup D_n^{d_0}\subseteq D$, for every $d_0\in \overline{D}_{p_1}$. This further implies that 
$\bigcup_{d_0\in \overline{D}_{p_1}} D_n^{d_0}=D_n\subseteq D$, which is a contradiction.



\smallskip

Now,  by employing the reasoning used in  Case 2  of the proof of Lemma \ref{lem:mineig}, we can see that the
equality $\sum_{d\in D_n^{d_0}\cap D} c(j,\frac{n}{d}) =
-\frac{\varphi(\frac{n}{d_0})}{p_1-1}$ holds for every $d_0\in \overline{D}_{p_1}$ if and only if the conditions
$D_n^{d_0}\cap D=d_0$ ($\beta_1=0$), $\gamma_1=\alpha_1-1$ and $\frac{n_1}{d_0\gcd(n_1/d_0,j_1)}=1$ are
satisfied. 
We deduce that $D=\overline{D}_{p_1}$ from the facts that $\overline{D}_{p_1}\subseteq D$ and $D_n^{d_0}\cap D=d_0$, for all $d_0\in \overline{D}_{p_1}$.

Moreover, from $\frac{n_1}{d_0\gcd(n_1/d_0,j_1)}=1$ we have
$\frac{n_1}{d_0}\mid j_1$ for all $d_0\in \overline{D}_{p_1}$. Finally,
since $\lcm\{\frac{n_1}{d_0}\ |\ d_0\in \overline{D}_{p_1} \}=n_1$ we see that $n_1\mid j_1$. Additionally, as $p_1^{\alpha_1-1}\| j$ ($\gamma_1=\alpha_1-1$)  we
conclude that the equality holds for every $j$ such that
$p_1^{\alpha_1-1}n_1\mid j$, and thus $\frac{n}{p_1}\mid j$.
\end{dok}

 By the previous assertion we determine the unique graph, which is $\ICG_n(D)$ and $D=\{d\in D_n\ |\ p_1\nmid d\}$. This particular graph possesses the lowest eigenvalue when compared to all connected integral circulant graphs of order $n$. This eigenvalue, represented as $\lambda_j(n,D)$, is precisely equal to $-\frac{n}{p_1}$ and can be derived for any index $j$ that is a multiple of $\frac{n}{p_1}$ within the sequence of eigenvalues as defined in (\ref{ldef}).
 Furthermore, the index of the this graph can be calculated in the following way
$$
\sum_{p_1\nmid d\atop d \in D_n}\varphi(\frac{n}{d})=\sum_{d \in D_n}\varphi(\frac{n}{d})-\sum_{p_1\mid d\atop d \in D_n}\varphi(\frac{n}{d})=n-\sum_{d' \in D_{\frac{n}{p_1}}}\varphi(\frac{n}{p_1d'})=n-\frac{n}{p_1}. 
$$
If we denote the spread of the integral circulant graph $\ICG_n(D)$ by $s(\ICG_n(D))\\=\max\{\lambda_0-\lambda_j\ |\ 1\leq j\leq n-1\}$ we see that $s(\ICG_n(\{d\in D_n\ |\ p_1\nmid d\}))=n$.
We have already mentioned that the maximum spread problem (finding maximal spread in a class of connected graphs with a prescribed order) is closely related to the maximum index problem and the minimal least eigenvalue problem in a given class of graphs.

We will now proceed to calculate the maximal spread for all integral circulant graphs of a given order $n$, denoted as $\maxspread(n)$, where $\maxspread(n)$ is defined as the maximum spread among all integral circulant graphs of order $n$, that is, $\maxspread(n)=\max\{s(\ICG_n(D))\ |\ D\subseteq D_n\}$. Our aim is to demonstrate that $\maxspread(n)$ is indeed equal to $n$. Consequently, we can conclude that the graph denoted as $\ICG_n({d\in D_n\ |\ p_1\nmid d})$ belongs to the class of graphs that achieve this maximum spread.
In fact, we will provide a complete characterization of the integral circulant graphs $\ICG_n(D)$ that achieve this maximum spread among all integral circulant graphs of order $n$.
Furthermore, it is worth noting that the integral circulant graph of order $n$ with maximal index, which is the complete graph $K_n$, also has a spread equal to $n$.

\begin{thm}
For a positive integer $n>1$, it holds that
$\maxspread(n)=n$. The spread of
connected graph $\ICG_n(D)$ attains $\maxspread(n)=n$ if and only
if $\gcd(\{d\ |\ d\in D_n\setminus D\})>1$.
\end{thm}
\begin{dok}
Considering the well-established fact that the complete graph $K_n$ is isomorphic to $\ICG_n(D_n)$, as indicated by (\ref{ldef}), it becomes evident that $\lambda_j(n,D)=\sum_{d\in D_n}c(j,\frac{n}{d})=-1$, for any value of $j$ within the range of $1\leq j\leq n-1$.
Moreover, when we take into account that $\sum_{d\in D_n}c(j,\frac{n}{d})=-1$ holds true, then for any graph $\ICG_n(D)$ and for any integer $j$ satisfying $1\leq j\leq n-1$, it becomes apparent that  
\begin{eqnarray*}
\lambda_0-\lambda_j&=&\sum_{d\in D}\varphi(\frac{n}{d})-\sum_{d\in D}c(j,\frac{n}{d})\\
&=&\sum_{d\in D}\varphi(\frac{n}{d})+1+\sum_{d\in D_n\setminus D}c(j,\frac{n}{d})\leq \sum_{d\in
D_n}\varphi(\frac{n}{d})+1=n.
\end{eqnarray*}
Now we find all the divisor sets $D$ such that $s(\ICG_n(D))=n$, for a given order $n$.
First, we prove the following assertion. Let $\ICG_n(D)$ be an
integral circulant graph such that $D=\{d_1,\ldots,d_t\}$ and
$\gcd(d_1,\ldots,d_t)=d$. Then $d>1$ if and only if
$\lambda_0=\lambda_{\frac{n}{d}}$.

Let $\ICG_n(D)$ be an integral circulant graph such that
$\gcd(d_1,\ldots,d_t)=d>1$, where $D=\{d_1,\ldots,d_t\}$. For $j=\frac{n}{d}$, we have that $\frac{n}{d_i}\mid \frac{n}{d}$,  $t_{\frac
n{d_i},j}=\frac{n}{d_i\gcd(\frac{n}{d_i},j)}=1$, resulting in $c(j,\frac{n}{d_i})=\varphi(\frac{n}{d_i})$, for every $1\leq i\leq t$. Consequently, 
$\lambda_{\frac{n}{d}}=\sum_{i=1}^t\varphi(\frac{n}{d_i})=\lambda_0$.
On the other hand, if $\lambda_0=\lambda_j$ then $\sum_{i=1}^t
\varphi(\frac{n}{d_i})-c(j,\frac{n}{d_i})=0$. Moreover, since every
summand $\varphi(\frac{n}{d_i})-c(j,\frac{n}{d_i})$ is greater or
equal to zero, we have that
$\varphi(\frac{n}{d_i})=c(j,\frac{n}{d_i})$, for every $d_i\in D$.
This is equivalent to $\mu(t_{\frac{n}{d_i},j})=1$ and $\varphi(t_{\frac{n}{d_i},j})=1$ and hence to $t_{\frac{n}{d_i},j}=1$.
Considering the definition of $t_{\frac{n}{d_i},j}$ in (\ref{ramanujan}), we observe that $t_{\frac{n}{d_i},j}=1$ if and only if $\gcd(\frac{n}{d_i},j)=\frac{n}{d_i}$. This condition is satisfied when $\frac{n}{d_i}\mid j$ for $1\leq i\leq t$. Consequently, the previous equality holds if and only if $\lcm(\frac{n}{d_1},\ldots,\frac{n}{d_t})\mid j$. On the
other hand, if $\gcd(d_1,\ldots,d_t)=1$ then
$\lcm(\frac{n}{d_1},\ldots,\frac{n}{d_t})=n\mid j$, which is a
contradiction since the index $j$ is supposed to take values from the range
$1\leq j\leq n-1$. Therefore, $\gcd(d_1,\ldots,d_t)>1$.

Now let $\lambda_0,\ldots,\lambda_{n-1}$ and
$\mu_0,\ldots,\mu_{n-1}$ be the eigenvalues of $\ICG_n(D)$ and
$\ICG_n(D_n\setminus D)$, respectively. Since $\ICG_n(D_n\setminus
D)$ is the complement of $\ICG_n(D)$ it is well known that
$\mu_0=n-1-\lambda_0$ and $\mu_i=-1-\lambda_i$, for $1\leq i\leq
n-1$. Based on the discussion in the previous paragraph, when considering $\ICG_n(D_n\setminus D)$, we can conclude that $d=\gcd\{d\ |\ d\in D_n\setminus D\}>1$ if and only if $\mu_0=\mu_{\frac{n}{d}}$. This is equivalent to
$n-1-\lambda_0=-1-\lambda_{\frac{n}{d}}$, which implies that
$\lambda_0-\lambda_{\frac{n}{d}}=n\leq s(\ICG_n(D))\leq
\maxspread(n)=n$. The last chain of inequalities evidently shows that $s(\ICG_n(D))=n$.
\end{dok}

\medskip

In the following assertion, we will identify all graphs in which the least eigenvalue is the smallest among all connected integral circulant graphs of a given order, while provided that that their complements are also connected.
Throughout the remainder of this section, we will continue to use the same notation as established in Lemma \ref{lem:mineig}: $S_{p_1}(n)=\alpha_1$, $S_{p_1}(j)=\gamma_1$, $S_{p_1}(d)=\beta_1$, $n_1=\frac{n}{p_1^{\alpha_1}}$,
$j_1=\frac{j}{p_1^{\gamma_1}}$, $d_0=\frac{d}{p_1^{\beta_1}}$, $N_{d_0}=N(n_1, d_0, j_1)=\frac{n_1}{d_0\gcd(j_1,\frac{n_1}{d_0})}$ and $N=N(n_1, d_0', j_1)=\frac{n_1}{d_0'\gcd(j_1,\frac{n_1}{d_0'})}$.

\begin{thm}
\label{minimal_least_complement}
Let $\ICG_n(D)$ be an arbitrary connected integral circulant graph of order $n$ whose complement is also connected. Then the
following inequality holds for every $1\leq j\leq n$
and $D\subseteq D_n$
$$
\lambda_j(n,D)\geq -\frac{n} {p_1}+p_1^{\alpha_1-1}.
$$ 
The equality is satisfied in the following cases

\begin{itemize}
\item[i)] $n$ is an arbitrary positive integer,  $D=\overline{D}_{p_1}\setminus \{\frac{n}{p_1^{\alpha_1}}\}$ and $\frac {n}{p_1}\mid j$, for $1\leq j\leq
n-1$,
\item [ii)] $n=2^{\alpha_1}3$, $D=\{2,4,\ldots,2^{\alpha_1},3\}$, $2^{\alpha_1-1}\parallel j$ and $3\nmid j$, for $1\leq j\leq
n-1$.
\end{itemize}

\label{thm:both_connected}

\end{thm}

\begin{dok}
According to the statement of Theorem \ref{thm:mleigenvalues}, for $\ICG_n(D)$ where $D=\overline{D}_{p_1}\setminus \{\frac{n}{p_1^{\alpha_1}}\}$, and $1\leq j\leq n-1$ with $\frac {n}{p_1}\mid j$, the eigenvalues are given by
\begin{equation}
\label{equality}
\lambda_j(n,D)=\sum_{d\in \overline{D}_{p_1}} c(j,\frac{n}{d})-c(j,p_1^{\alpha_1})=-\frac{n}{p_1}+p_1^{\alpha_1-1}.
\end{equation}

Indeed, we have that $t_{p_1^{\alpha_1},j}=\frac{p_1^{\alpha_1}}{\gcd(p_1^{\alpha_1}, j)}=p_1$ and $c(j,p_1^{\alpha_1})=-\frac{\varphi(p_1^{\alpha_1})}{\varphi(p_1)}=-p_1^{\alpha_1-1}$.
It is evident that $\ICG_n(D_n\setminus D)$ is a connected graph since $p_1$ and $\frac{n}{p_1^{\alpha_1}}$ belong to $D_n\setminus D$, ensuring that $\gcd(\{d\in D_n\setminus D\})=1$.

In the rest of the proof, we will show that $\lambda_j(n,D)\geq -\frac{n}{p_1}+p_1^{\alpha_1-1}$, for $0\leq j\leq n-1$ and any $D\subseteq D_n$ such that $\gcd(\{d\in D_n\setminus D\})=1$.
Furthermore, we will determine necessary conditions for which the equality is achieved.
We distinguish three cases depending on the value of $\gamma_1$.

{\noindent \bf Case 1. } $\gamma_1\leq \alpha_1-2$.
We will prove that $\lambda_j(n,D)>  -\frac{n}{p_1}+p_1^{\alpha_1-1}$, for any $D\subseteq D_n$, $0\leq j\leq n-1$ and $\gamma_1\leq \alpha_1-2$.

We may notice that $\lambda_j(n,D_0)$ attains minimal value either in Case 2 or Case 3 from the proof of Lemma \ref{lem:mineig}, for $D_0\subseteq D_n^{d_0}$ and some $d_0\in \overline{D}_{p_1}$.  In either scenario, utilizing either relation (\ref{eq:sum c}) or (\ref{eq:sumc1}), we directly conclude that 
\begin{equation}
 \sum_{d\in D_0} c(j,\frac{n}{d})\geq-\frac{\varphi(\frac{n_1}{d_0})}{\varphi(N_{d_0})}p_1^{\gamma_1}.
 \label{eq:for c}
\end{equation}
Furthermore, as $\gamma_1\leq \alpha_1-2$ we see that $\sum_{d\in D_0} c(j,\frac{n}{d})\geq-\frac{\varphi(\frac{n_1}{d_0})}{\varphi(N_{d_0})}p_1^{\alpha_1-2}$ and therefore 
\begin{eqnarray*}
\lambda_j(n,D)&=&
\sum_{d_0\in\overline{D}_{p_1}\cap D} \sum_{d\in D_n^{d_0}\cap D}
c(j,\frac{n}{d})\geq -\sum_{d_0\in\overline{D}_{p_1}}
\frac{\varphi(\frac{n_1}{d_0})}{\varphi(N_{d_0})}p_1^{\alpha_1-2}\\
&=&-p_1^{\alpha_1-2}\sum_{d_0\in\overline{D}_{p_1}}
\varphi(\frac{n_1}{d_0})
=-p_1^{\alpha_1-2}n_1=-\frac
{n} {p_1^2},
\end{eqnarray*}
for some $D\subseteq D_n$. Now, given that $-\frac{n}{p_1^2}>-\frac{n}{p_1}+p_1^{\alpha_1-1}$ is equivalent to $n(p_1-1)>p_1^{\alpha_1+1}$, we obtain that $\lambda_j(n,D)>-\frac{n}{p_1}+p_1^{\alpha_1-1}$ holds true when $n\neq p_1^{\alpha_1}$. However, when $n = p_1^{\alpha_1}$, it follows that either $\ICG_n(D)$ or $\ICG_n(D_n\setminus D)$ must be disconnected. This is due to the fact that if $1\in D_n\setminus D$, then $p_1\mid d$ for every $d\in D$, and conversely, if $1\in D$, then $p_1\mid d$ for every $d\in D_n\setminus D$. As a result, this scenario is excluded from our consideration.

{\noindent \bf Case 2. } $\gamma_1=\alpha_1$. 
Since $\gamma_1=\alpha_1$, it follows that the inequality $\beta_1\geq \alpha_1-\gamma_1$ trivially holds for all $d\in D_0\subseteq D_n^{d_0}$ and $d_0\in \overline{D}_{p_1}$, and thus we can deduce that $\sum_{d\in D_0} c(j,\frac{n}{d})\geq-\frac{\varphi(\frac{n_1}{d_0})}{\varphi(N_{d_0})}p_1^{\alpha_1}$,  as indicated by Case 3 of Lemma \ref{lem:mineig}. However, the last inequality can not be achieved for $d_0=\frac{n}{p_1^{\alpha_1}}=n_1$ due to $N_{d_0}=t_{\frac{n_1}{d_0},j_1}=\frac{1}{\gcd(1,j_1)}=1$ and $\mu(N_{d_0})=1$. As a result, we can deduce that $c(j,\frac{n}{d})\geq 0$ holds true for $d\in D_n^{n/p_1^{\alpha_1}}$. Based on the preceding analysis, for $D\subseteq D_n$ and $\frac{n}{p_1^{\alpha_1}}\not\in D$, it holds that

\begin{eqnarray*}
\lambda_j(n,D)&=&
\sum_{d_0\in\overline{D}_{p_1}\setminus \{\frac{n}{p_1^{\alpha_1}}\}\cap D} \sum_{d\in D_n^{d_0}\cap D}
c(j,\frac{n}{d})\geq -\sum_{d_0\in\overline{D}_{p_1}\setminus \{\frac{n}{p_1^{\alpha_1}}\}}
\frac{\varphi(\frac{n_1}{d_0})}{\varphi(N_{d_0})}p_1^{\alpha_1}.\\
\end{eqnarray*}

As $p_1\leq p_2-1\leq \varphi(N_{d_0})$, we see that 
\begin{eqnarray*}-\sum_{d_0\in\overline{D}_{p_1}\setminus \{\frac{n}{p_1^{\alpha_1}}\}}
\frac{\varphi(\frac{n_1}{d_0})}{\varphi(N_{d_0})}p_1^{\alpha_1}&\geq& -\sum_{d_0\in\overline{D}_{p_1}\setminus \{\frac{n}{p_1^{\alpha_1}}\}}
\varphi(\frac{n_1}{d_0})p_1^{\alpha_1-1}\\
&=& -(n_1-1)p_1^{\alpha_1-1}=-\frac{n}{p_1}+p_1^{\alpha_1-1}.
\end{eqnarray*}

This means that 
$\lambda_j(n,D)\geq -\frac{n}{p_1}+p_1^{\alpha_1-1}$ and the equality holds if and only if $D=\cup_{d_0\in \overline{D}_{p_1}\setminus \{\frac{n}{p_1^{\alpha_1}}\}}  D_n^{d_0}$, $N_{d_0}=3$ and $p_1=2$. However, under those conditions, we find that  $D_n\setminus D=D_{n}^{n/p_1^{\alpha_1}}$,  leading to the conclusion that  $\ICG_n(D_n\setminus D)$ is disconnected, as $\gcd(\{d\in D_n^{n/p_1^{\alpha_1}}\})=\frac{n}{p_1^{\alpha_1}}$, which is a contradiction with the statement of the assertion.

{\noindent \bf Case 3. } $\gamma_1=\alpha_1-1$. According to (\ref{eq:for c}), we see that $\sum_{d\in D_0} c(j,\frac{n}{d})\geq-\frac{\varphi(\frac{n_1}{d_0})}{\varphi(N_{d_0})}p_1^{\alpha_1-1}$,
for $D_0\subseteq D_n^{d_0}$ and some $d_0\in \overline{D}_{p_1}$.
Notice fist that there exists a $d'_0\in \overline{D}_{p_1}$ such that $d'_0\not \in D$.  If the opposite were true, all divisors from $D_n\setminus D$ would be divisible by $p_1$, which contradicts the assumption that $\ICG_n(D_n\setminus D)$ is connected.

Suppose  that $D_n^{d'_0}\cap D= \emptyset$. 

The subsequent sequence of inequalities is valid for any $1\leq j\leq n-1$
\begin{eqnarray}
\label{eq:first}
&&
\lambda_j(n,D)=\sum_{d_0\in(\overline{D}_{p_1}\setminus \{d'_0\})\cap D} \sum_{d\in D_n^{d_0}\cap D}
c(j,\frac{n}{d})\geq -\sum_{d_0\in\overline{D}_{p_1}\setminus \{d_0'\}}
\frac{\varphi(\frac{n_1}{d_0})}{\varphi(N_{d_0})}p_1^{\alpha_1-1} \nonumber\\
&\geq&-\sum_{d_0\in\overline{D}_{p_1}}
\varphi(\frac{n_1}{d_0})p_1^{\alpha_1-1}+\varphi(\frac{n_1}{d_0'})p_1^{\alpha_1-1}=
-\frac{n}{p_1}+\varphi(\frac{n_1}{d_0'})p_1^{\alpha_1-1}.
\end{eqnarray}

Due to the fact that $\lambda_j(n,D)=\sum_{d_0\in(\overline{D}_{p_1}\setminus \{d'_0\})\cap D} \sum_{d\in D_n^{d_0}\cap D}
c(j,\frac{n}{d})$, we have that $\lambda_j(n,D)\geq -\frac{n}{p_1}+\varphi(\frac{n_1}{d_0'})p_1^{\alpha_1-1}$ and the right hand side of the equality is always greater than or equal to $-\frac{n}{p_1}+p_1^{\alpha_1-1}$. The equality holds if and only if $\varphi(\frac{n_1}{d_0'})=1$, and 
for such $j$ and $D\cup \{d_0'\}$ that satisfy  the equality in (\ref{ineq:lambdaj}). Since $n_1$ is an odd number, it follows that $\frac{n_1}{d_0'}=1$, which in turn implies $d_0'=n_1=\frac{n}{p_1^{\alpha_1}}$.
Furthermore, based on the proof of Theorem \ref{thm:mleigenvalues}, the condition of equality in (\ref{ineq:lambdaj}) is satisfied exclusively when $D=\overline{D}_{p_1}\setminus \{\frac{n}{p_1^{\alpha_1}}\}$ and $\lcm\{\frac{n_1}{d_0'}\ |\ d_0'\in \overline{D}_{p_1}\}=n_1\mid j_1$. This conclusion concludes the fist part of the proof.

\medskip

Suppose now that $D_n^{d'_0}\cap D\neq \emptyset$.

Once again, according to (\ref{eq:for c}),  we have that  
\begin{eqnarray}
\label{eq:second}
\sum_{d\in D_n^{d'_0}\cap D}
c(j,\frac{n}{d})\geq -\frac{\varphi(\frac{n_1}{d'_0})}{\varphi(N)}p_1^{\alpha_1-1}.
\end{eqnarray}

Recall that we denoted $N$ to be $N(n_1,d'_0,j_1)$ and $N_{d_0}$ be $N(n_1,d_0,j_1)$ for $d_0\neq d_0'$ and $d_0\in \overline{D}_{p_1}$.
As $\lambda_j(n,D)=\sum_{d_0\in\overline{D}_{p_1}\cap D} \sum_{d\in D_n^{d_0}\cap D}
c(j,\frac{n}{d})$,  by summing the respective sides of inequalities  (\ref{eq:first}) and (\ref{eq:second})
we obtain that 
\begin{eqnarray}
\label{eq:ineq_lambda_j}    
\lambda_j(n,D)&\geq& -\frac{n}{p_1}+\varphi(\frac{n_1}{d_0'})p_1^{\alpha_1-1}-\frac{\varphi(\frac{n_1}{d'_0})}{\varphi(N)}p_1^{\alpha_1-1}\nonumber\\
&=&-\frac{n}{p_1}+p_1^{\alpha_1-1}\varphi(\frac{n_1}{d'_0})\frac{\varphi(N)-1}{\varphi(N)}.
\end{eqnarray}

Moreover, since $\varphi(\frac{n_1}{d'_0})\frac{\varphi(N)-1}{\varphi(N)}\geq 1$ ($\varphi(\frac{n_1}{d'_0})\geq 2\geq \frac{\varphi(N)}{\varphi(N)-1}$) we have that $\lambda_j(n,D)\geq -\frac{n}{p_1}+p_1^{\alpha_1-1}$. The equality holds if and only if 
$\varphi(\frac{n_1}{d'_0})(\varphi(N)-1)=\varphi(N)$, $N_{d_0}=1$, $D\cap (\cup_{d_0\neq d_0'} D_n^{d_0})=\overline{D}_{p_1}\setminus\{d'_0\}$ and $D\cap  D_n^{d'_0}=D_n^{d'_0}\setminus \{d'_0\}$.
However, as $\gcd(\varphi(N)-1,\varphi(N))=1$, we have that the last equality is satisfied if and only if $\varphi(N)=\varphi(\frac{n_1}{d'_0})=2$. As $2=\varphi(N)\geq p_2-1\geq p_1$, we have that $p_1=2$ and $p_2=N=3$. Moreover, since $n_1/d'_0$ is odd we have that $d'_0=\frac{n_1}{3}$.
The equality $N=3$ implies that 
$t_{\frac{n_1}{d_0'},j_1}=\frac{n_1}{d_0'\gcd(\frac{n_1}{d_0'},j_1)}=\frac{3}{\gcd(3,j_1)}=3$ or equivalently $\gcd(3,j_1)=1$ which means $3\nmid j_1$.
On the other hand, from $N_{d_0}=\frac{n_1}{d_0\gcd(\frac{n_1}{d_0},j_1)}=1$, we see that $\frac{n_1}{d_0}\mid j_1$.
If $d_0'\neq 1$, then we obtain that $N_1=1$ implying that $n_1\mid j_1$. This conclusion contradicts the fact that $3\nmid j_1$, which means that $d_0'=\frac{n_1}{3}=1$.
Since $n_1=3$ we obtain that $d_0=3$ is the only divisor distinct from $d_0'$ such that $d_0\in \overline{D}_{p_1}$.
Therefore, the equality $N_{d_0}=\frac{1}{\gcd(1,j_1)}=1$ is always satisfied, which means that the second part of the assertion is finished. 


\end{dok}

As per the initial part of Theorem \ref{thm:both_connected}, it becomes evident that within the class of integral circulant graphs $\ICG_n(D)$ having order $n=2^{\alpha_1}3$, where both graphs $\ICG_n(D)$ and $\overline{\ICG_n(D)}$ are connected, the spectrum of the unitary Cayley graph $\ICG_n({1})$ includes the minimal least eigenvalue. This holds true since $\overline{D}_{p_1}=\{1,3\}$ and $\frac{n}{p_1^{\alpha_1}}=3$.
The second part of the theorem presents an additional graph of order $n=2^{\alpha_1}3$, as expected.


\medskip

In the subsequent statement, we will identify the second minimal least eigenvalue within the spectra of all connected integral circulant graphs with an order of $n$ and furthermore, we will determine all the graphs in this class that possess the second minimal least eigenvalue. Given that $\ICG_n(\overline{D}_{p_1})$ represents the unique graph with the minimal eigenvalue among its spectra, the problem of finding the second minimal least eigenvalue can be reduced to determining the minimal least eigenvalue within the class of all connected integral circulant graphs $\ICG_n(D)$ of a given order $n$ with sets $D$ distinct from $\overline{D}_{p_1}$.
Furthermore, when $n$ is a prime number, there exists a unique connected integral circulant graph denoted as $\ICG_n(1)$. Consequently, according to the definition of the second minimal least eigenvalue, the concept of a second minimal least eigenvalue does not apply in this context. Therefore, in the following theorem we assume that $n$ is a composite number.

\begin{thm}
Let $\ICG_n(D)$ be any  connected integral circulant graph of the order $n$, where $n$ is a composite number and $D\neq \overline{D}_{p_1}$. Then the
following inequality holds for every $1\leq j\leq n$
and $D\subseteq D_n$
\begin{eqnarray*}
\lambda_j(n,D) &\geq&\left\{
\begin{array}{rl}
-\frac{n} {p_1}+p_1-1, &  \alpha_1> 1\\
-\frac{n} {p_1}+1, &  \alpha_1= 1 \\
\end{array}\right..
\end{eqnarray*}

Equality is achieved under the following conditions

\begin{itemize}
\item[i)] when $\alpha_1>1$,  $D=\overline{D}_{p_1}\cup \{\frac{n}{p_1}\}$ and $\frac {n}{p_1}\mid j$, for $1\leq j\leq
n-1$,
\item[ii)] when $\alpha_1=1$,  $D=\overline{D}_{p_1}\setminus \{\frac{n}{p_1}\}$ and $\frac {n}{p_1}\mid j$, for $1\leq j\leq
n-1$,
\item [iii)] for $n=6$, when $D=\{1\}$ and $j=3$, or $D=\{1,2\}$ and $j=\{2,4\}$, or $D=\{2,3\}$ and $j\in\{1,5\}$.
\end{itemize}

\label{thm:second_smallest}

\end{thm}

\begin{dok}
To identify the second minimal least eigenvalue within the integral circulant graphs of order $n$, one can simply find the minimal least eigenvalue among the integral circulant graphs denoted as $\ICG_n(D)$, where $D\neq \overline{D}_{p_1}$ for a given order $n$.
This conclusion is based on the fact that $\ICG_n(\overline{D}_{p_1})$ is the only graph whose spectrum includes the minimal least eigenvalue.

In accordance with (\ref{equality}), it becomes evident that the eigenvalue $\lambda_j$ present in the spectrum of the graph $\ICG_n(\overline{D}_{p_1}\setminus \{\frac{n}{p_1^{\alpha_1}}\})$, for $1\leq j\leq n-1$ where $\frac {n}{p_1}\mid j$, equals $-\frac{n}{p_1}+1$, in the case where $\alpha_1=1$. On the other hand, by employing a similar analysis to what was conducted in equation (\ref{equality}), we can demonstrate that the eigenvalue $\lambda_j$ within the spectrum of the graph $\ICG_n(\overline{D}_{p_1}\cup \{\frac{n}{p_1}\})$, where $1\leq j\leq n-1$ and $\frac {n}{p_1}\mid j$, is equivalent to $-\frac{n}{p_1}+c(j,p_1)$, under the condition that $\alpha_1>1$.
Moreover, as  $t_{p_1,j}=\frac{p_1}{\gcd(p_1, j)}=1$ and $c(j,p_1)=\varphi(p_1)=p_1-1$, we conclude that $\lambda_j=-\frac{n}{p_1}+p_1-1$, for $\alpha_1>1$.

In the remainder of the proof, we will demonstrate the inequality stated in the assertion for any $D\neq \overline{D}_{p_1}$ and $0\leq j\leq n-1$.
Furthermore, we will establish the necessary conditions under which equality is attained.
We distinguish three cases depending on the value of $\gamma_1$ as previously demonstrated in the preceding theorem.

{\noindent \bf Case 1. } $\gamma_1\leq \alpha_1-2$.
In Case 1 of the previous theorem's proof, it has been demonstrated that  $\lambda_j(n,D)\geq  -\frac{n}{p_1^2}$, for any $D\subseteq D_n$, $0\leq j\leq n-1$ and $\gamma_1\leq \alpha_1-2$.
Now, starting from the inequality $-\frac{n}{p_1^2}>-\frac{n}{p_1}+p_1-1$, which can be rewritten as $n>p_1$, we can deduce that $\lambda_j(n,D)>-\frac{n}{p_1}+p_1-1\geq -\frac{n}{p_1}+1$, given that $n$ is a composite number.

{\noindent \bf Case 2. } $\gamma_1=\alpha_1$.
In accordance with Case 2 of the proof of Theorem \ref{minimal_least_complement}, we conclude that 
$\lambda_j(n,D)\geq  -\frac{n}{p_1}+p_1^{\alpha_1-1}$, for any $D\subseteq D_n$, $0\leq j\leq n-1$ and $\gamma_1= \alpha_1$.
Moreover, given that $\alpha_1 > 1$, we can observe that $-\frac{n}{p_1} + p_1^{\alpha_1-1} > -\frac{n}{p_1} + p_1 - 1$. Consequently, it follows that $\lambda_j > -\frac{n}{p_1} + p_1 - 1$ for all $D\subseteq D_n$,
$0\leq j\leq n-1$, $\gamma_1= \alpha_1$ and $\alpha_1>1$. 

However, when $\alpha_1=1$, we can observe that $\lambda_j(n,D)\geq -\frac{n}{p_1}+1$. Consequently, we are able to identify the necessary conditions under which equality is achieved in the preceding inequality.
According to the proof outlined in Theorem \ref{minimal_least_complement} (Case 2), equality is achieved under the following conditions: when $p_1=2$,  $D=\cup_{d_0\in \overline{D}_{p_1}\setminus \{\frac{n}{p_1^{\alpha_1}}\}}  D_n^{d_0}$ and $N_{d_0}=3$ for all $d_0\in  \overline{D}_{p_1}\setminus \{\frac{n}{p_1^{\alpha_1}}\}$.
As $N_{d_0}=3$, it becomes evident that $p_2=3$.
Furthermore, we notice that $1\in \overline{D}_{p_1}\setminus \{\frac{n}{p_1^{\alpha_1}}\}$, as otherwise, it must be $1=\frac{n}{p_1}$, which would result in a contradiction since $n$ is a composite number.
Hence, it follows that $N_1 = \frac{n_1}{\gcd(j_1, n_1)} = 3$, which leads us to the conclusion that $\frac{n_1}{3}\mid j_1$.
Given that $\frac{n_1}{3}\mid j_1$, if we make the assumption that $p_2$ belongs to set $\overline{D}_{p_1}\setminus \{\frac{n}{p_1^{\alpha_1}}\}$, we can deduce that $N_{p_2} = \frac{n_1}{3\gcd(j_1, \frac{n_1}{3})} = 1$, which is a contradiction since our initial assumption was that  $N_{p_2}=3$. Hence, from $p_2\not\in \overline{D}_{p_1}\setminus \{\frac{n}{p_1^{\alpha_1}}\}$, we conclude that $p_2=\frac{n}{p_1}$, which is equivalent to $n=6$. Since $\overline{D}_{p_1}=\{1,3\}$, we have that $D=\cup_{d_0\in \overline{D}_{p_1}\setminus \{\frac{n}{p_1^{\alpha_1}}\}}  D_n^{d_0}=D_{6}^{1}=\{1,2\}$. Finally, given the fact that $\gamma_1=\alpha_1=1$ and $1=\frac{n_1}{3}\mid j_1$, we can conclude that $j$ must be in the set $\{2,4\}$ as these are the values of $j$ for which $\lambda_j$ achieves the value of $-\frac{n}{p_1}+1=-2$. In this manner, we ascertain the condition described in the second part of item $(iii)$ of the assertion, under which the inequality is satisfied.

{\noindent \bf Case 3. } $\gamma_1=\alpha_1-1$.

Suppose now $D\subseteq \overline{D}_{p_1}$. This means that there exists $d_0'\in \overline{D}_{p_1}$ such that $d_0'\not\in D$, as $\overline{D}_{p_1}\neq D$.
Given that $0=\beta_1=\alpha_1-\gamma_1-1$, for all $d\in D$, according to Case 2 of Lemma \ref{lem:mineig}, we conclude that 
$$
\lambda_{j}(n,D)\geq -\sum_{d\in D\setminus\{d_0'\}}p_1^{\alpha_1-1}\varphi(\frac{n_1}{d})=-\frac{n}{p_1}+p_1^{\alpha_1-1}\varphi({\frac{n_1}{d'_0}})\geq -\frac{n}{p_1}+p_1^{\alpha_1-1}.
$$
Consequently, it follows that $\lambda_{j}(n,D) > -\frac{n}{p_1} + p_1 - 1$ for all $D\subseteq D_n$,
$0\leq j\leq n-1$, $\gamma_1= \alpha_1-1$ and $\alpha_1>1$.
However, when $\alpha_1=1$, we observe  that  $\lambda_{j}(n,D)\geq -\frac{n}{p_1}+1$, and the equality holds for $n_1=d_0'$, implying that $d_0'=\frac{n}{p_1^{\alpha_1}}=\frac{n}{p_1}$.
In this way, we demonstrate the validity of item $(ii)$ and the first part of item $(iii)$ in the theorem's statement. 

Suppose now that $D\not\subseteq \overline{D}_{p_1}$. This implies that there exists $d_0\in D$ such that $d_0\not\in \overline{D}_{p_1}$.
Let $d_0'$ be a divisor such that $d_0\in D_{n}^{d_0'}$. 

If $d_0'\not\in D$, then, based on (\ref{eq:ineq_lambda_j}),  we can establish that $\lambda_j(n,D)\geq -\frac{n}{p_1}+p_1^{\alpha_1-1}\varphi(\frac{n_1}{d'_0})\frac{\varphi(N)-1}{\varphi(N)}\geq -\frac{n}{p_1}+p_1^{\alpha_1-1}$. Moreover, when $\alpha_1>1$, we observe that $\lambda_j(n,D)\geq -\frac{n}{p_1}+p_1^{\alpha_1-1}\geq -\frac{n}{p_1}+p_1>-\frac{n}{p_1}+p_1-1$.
In the case of  $\alpha_1=1$, we conclude that $\lambda_j(n,D)\geq -\frac{n}{p_1}+1$, which  holds true when $\varphi(\frac{n_1}{d'_0})\frac{\varphi(N)-1}{\varphi(N)}=1$. 
After examining the paragraph following (\ref{eq:ineq_lambda_j}), we find that the equation $\lambda_j(n,D) = -\frac{n}{p_1} + 1$ is satisfied under the conditions outlined in item $(ii)$ of Theorem \ref{minimal_least_complement}. 
This implies that the equality holds for $n=6$, $D=\{2,3\}$, and $j=\{1,5\}$, thereby satisfying the condition outlined in the third part of item $(iii)$ of the assertion.

If $d_0'\in D$, then, according to Cases 2 and 3 of Lemma \ref{lem:mineig}, we see that $c(j,\frac{n}{d_0'})$ and $c(j,\frac{n}{d_0})$ have the opposite signs ($S_{p_1}(d_0')=\alpha_1-\gamma_1-1=0$ and $S_{p_1}(d_0)>0=\alpha_1-\gamma_1-1$). 
Based on (\ref{eq:ineq_lambda_j}), we can notice that
\begin{eqnarray*}
\label{eq:inequality_}
\lambda_j(n,D)&\geq&  -\frac{n}{p_1}+\varphi(\frac{n_1}{d_0'})p_1^{\alpha_1-1}-\frac{\varphi(\frac{n_1}{d'_0})}{\varphi(N)}p_1^{\alpha_1-1}+c(j,\frac{n}{d_0'})\\
&=&-\frac{n}{p_1}+p_1^{\alpha_1-1}\varphi(\frac{n_1}{d'_0})\frac{\varphi(N)-1}{\varphi(N)}+c(j,\frac{n}{d_0'}).
\end{eqnarray*}

If $c(j,\frac{n}{d_0'})>0$, then we observe that $\lambda_j(n,D)>-\frac{n}{p_1}+p_1^{\alpha_1-1}> -\frac{n}{p_1}+p_1-1$ when $\alpha_1>1$, and $\lambda_j(n,D)>-\frac{n}{p_1}+1$ when $\alpha_1=1$.
If $c(j,\frac{n}{d_0'})=0$, then $\mu(N)=0$, and consequently, $c(j,\frac{n}{d})=0$ for every $d\in D_{n}^{d_0'}$. This, in turn, implies that $\sum_{d\in D_n^{d'_0}\cap D} c(j,\frac{n}{d})=0$, and as a result, we can enhance the lower bound for $\lambda_j(n,D)$ in the following manner: $\lambda_j(n,D)\geq -\frac{n}{p_1}+\varphi(\frac{n_1}{d_0'})p_1^{\alpha_1-1}$ (according to (\ref{eq:first})). For $\alpha_1>1$, it is evident that $\lambda_j(n,D) > -\frac{n}{p_1} + p_1 - 1$. In the case where $\alpha_1=1$, we have $\lambda_j(n,D) \geq -\frac{n}{p_1} + 1$, and equality is achieved when $\varphi(\frac{n_1}{d_0'}) = 1$. Consequently, we conclude that $d_0' = n_1$, which implies $N = \frac{1}{\gcd(1,j_1)} = 1$, leading to $\mu(N) = 1$. This contradicts the established fact that $\mu(N) = 0$, and hence, we must conclude that $\lambda_j(n,D) > -\frac{n}{p_1} + 1$.

Finally, if $c(j,\frac{n}{d_0'})<0$, we conclude that $c(j,\frac{n}{d_0})>0$.
Using (\ref{ineq:lambdaj}), we further have that $\lambda_j(n,D)\geq-\frac{n}{p_1}+c(j,\frac{n}{d_0})$. 
In accordance with equation (\ref{eq:Case 3}), we can derive the following lower bounds for $c(j, \frac{n}{d_0})$: $c(j, \frac{n}{d_0})\geq \min\{\frac{\varphi(\frac{n_1}{d_0'})}{\varphi(N)}, \frac{(p_1-1)\varphi(\frac{n_1}{d_0'})}{\varphi(N)}\}$, with these lower bounds being attained when $\beta_1=\alpha_1$ and $\beta_1=\alpha_1-1$, respectively, where $\beta_1=S_{p_1}(d_0)$. 
Given that $\frac{\varphi(\frac{n_1}{d_0'})}{\varphi(N)}\geq 1$ as per the definition of $N$, it becomes evident that $c(j, \frac{n}{d_0})\geq p_1-1$, provided that $\beta_1=\alpha_1-1$.
Furthermore, for $\beta_1=\alpha_1$, we can notice that $d_0'\neq n_1$. Otherwise, we would have that $d_0=p_1^{\beta_1}n_1=n$, and this divisor would  not belong to $D_n$, which leads to a contradiction. 
Considering that $N \mid \frac{n_1}{d_0'}$ and $\frac{n_1}{d_0'}$ is an odd number, it becomes evident that $\frac{\varphi(\frac{n_1}{d_0'})}{\varphi(N)} = 1$ if and only if $N = \frac{n_1}{d_0'}$. This, in turn, implies that $\gcd(\frac{n_1}{d_0'}, j_1) = 1$.
Referring to Theorem \ref{thm:mleigenvalues}, we can conclude that the lower bound in the inequality $\lambda_j(n, D) \geq -\frac{n}{p_1} + c(j, \frac{n}{d_0})$ is achievable if and only if $n_1 \mid j_1$.
Therefore, it follows that $\lambda_j(n, D)$ cannot equal $\frac{n}{p_1} + 1$, for $\beta_1=\alpha_1$, since the conditions $\gcd(\frac{n_1}{d_0'}, j_1) = 1$ and $n_1 \mid j_1$ cannot be simultaneously satisfied.

Based on the preceding discussion, we can draw the conclusion that when $\alpha_1=1$, the inequality simplifies to $\lambda_j(n, D) > -\frac{n}{p_1}+1$, since  $\beta_1>0$ holds (the case $\beta_1=\alpha_1$ must be satisfied).
However, for $\alpha_1>1$, we have that $c(j, \frac{n}{d_0})\geq \min\{\frac{\varphi(\frac{n_1}{d_0'})}{\varphi(N)}, p_1-1\}$, where $\frac{\varphi(\frac{n_1}{d_0'})}{\varphi(N)}> 1$.
Moreover, as $N\mid \frac{n_1}{d_0'}$ and $\frac{\varphi(\frac{n_1}{d_0'})}{\varphi(N)}> 1$, it can be readily obtained  that there exists a prime number $p_i$, where $2\leq i\leq k$, such that $p_i\mid \frac{\varphi(\frac{n_1}{d_0'})}{\varphi(N)}$ or $p_i-1\mid \frac{\varphi(\frac{n_1}{d_0'})}{\varphi(N)}$. This, in turn, implies that $c(j, \frac{n}{d_0})\geq \min\{\frac{\varphi(\frac{n_1}{d_0'})}{\varphi(N)}, p_1-1\}\geq \min\{p_i-1, p_1-1\}=p_1-1$. Therefore, we have established the following lower bound for $\lambda_j(n,D)$: $\lambda_j(n,D)\geq-\frac{n}{p_1}+p_1-1$, with equality holding if and only if both conditions $\gcd(\frac{n_1}{d_0'}, j_1) = 1$ and $n_1 \mid j_1$ are met. Both equation are satisfied if and only if $d_0'=n_1$, which further implies that $d_0=p_1^{\alpha_1-1}n_1=\frac{n}{p_1}$, given that $\beta_1=\alpha_1-1$. We have thus demonstrated the validity of assertion $(i)$ of the theorem in this way.
\end{dok}

\vspace{0.5cm}

\section*{Declarations}

\begin{itemize}
\item Funding

This research was supported by the research project of the Ministry of Education, Science and Technological Development of the Republic of Serbia (number 451-03-47/2023-01/ 200124).

\item Conflict of interest/Competing interests

The authors declare that they have no known competing financial interests or personal
relationships that could have appeared to influence the work reported in this paper.

\item Ethics approval

Not required.

\item Consent to participate

The requirement for informed consent was not necessary because the study relies on  theoretical concepts.

\item Consent for publication

This manuscript does not include details, images or videos relating to an individual person;
therefore, consent for publication is not required, beyond the informed consent provided by all
study participants as described above.

\item Availability of data and materials

Not Applicable

\item Code availability

Not Applicable

\item Authors' contributions

Conceptualization: MB. Methodology: MB. Formal analysis: MB. 
Writing — original draft preparation: MB. Writing — review and editing: MB.
This publication is the work of the author MB. All authors read and approved the final manuscript.

\end{itemize}

%
%



\end{document}